\documentclass[11pt]{article}%
\usepackage{amsmath}
\usepackage{amsfonts}
\usepackage{amssymb}
\usepackage{graphicx}
\usepackage{theorem}
\usepackage{makeidx}
\usepackage{verbatim}
\usepackage[colorlinks=true]{hyperref}%
\providecommand{\U}[1]{\protect\rule{.1in}{.1in}}
\newtheorem{thm}{Theorem}[section]

\newtheorem{pr}[thm]{Proposition}
\newtheorem{df}[thm]{Definition}

{\theorembodyfont{\upshape}

}
\numberwithin{equation}{section} 
\begin{document}

\title{ }

\begin{center}
\vspace*{1.5cm}

\textbf{ON PARAMETRIC\ VECTOR\ OPTIMIZATION VIA}

\textbf{METRIC\ REGULARITY OF CONSTRAINT\ SYSTEMS}

\vspace*{1cm}

M. DUREA

{\small {Faculty of Mathematics, "Al. I. Cuza" University,} }

{\small {Bd. Carol I, nr. 11, 700506 -- Ia\c{s}i, Romania,} }

{\small {e-mail: \texttt{durea@uaic.ro}}}

\bigskip

R. STRUGARIU

{\small Department of Mathematics, "Gh. Asachi" Technical University, }

{\small {Bd. Carol I, nr. 11, 700506 -- Ia\c{s}i, Romania,} }

{\small {e-mail: \texttt{rstrugariu@tuiasi.ro}}}
\end{center}

\bigskip

\bigskip

\noindent{\small {\textbf{Abstract: }}Some metric and graphical regularity
properties of generalized constraint systems are investigated. Then, these
properties are applied in order to penalize (in the sense of Clarke) various
scalar and vector optimization problems. This method allows us to present
several necessary optimality conditions in solid constrained vector
optimization.}

\bigskip

\noindent{\small {\textbf{Keywords: }}set-valued mappings $\cdot$ metric
regularity $\cdot$ scalar and vector optimization}

\bigskip

\noindent{\small {\textbf{Mathematics Subject Classification (2010): }90C30
$\cdot$ 49J53} {$\cdot$ 54C60}}

\section{Introduction and motivation}

In theory of optimization, a large literature is dedicated to the branch of
parametric optimization problems with equilibrium constraints and this is one
of the reasons which motivates the many efforts made in the last decades in
studying the variational systems (appearing as constraint systems in different
mathematical programs) in their own right and involving growing generality. We
quote here only the important works of Robinson \cite{Rob1980}, Dontchev and
Rockafellar \cite{DontRock2009b}, Mordukhovich \cite{Mor2006} for
comprehensive discussions and historical facts.

In this paper we firstly aim to underline some aspects of these topics by
presenting in a simple manner some optimization-related motivational facts
leading to two types of metric regularity for fairly general variational
systems. Then we survey, complete and extend some results previously obtained
by the authors. Lastly, we get optimality conditions for some vector
mathematical programs, by the use of the regularity of constraints system and
the power of Mordukhovich's generalized differentiation theory.

\bigskip

Let $X$ be a Banach space. In this setting, $B(x,r)$ and $D(x,r)$ denote the
open and the closed ball with center $x$ and radius $r,$ respectively$.$
Sometimes we write $\mathbb{D}_{X},\mathbb{B}_{X}$ and $\mathbb{S}_{X}$ for
the closed, the open unit ball and the unit sphere of $X,$ respectively. If
$x\in X$ and $A\subset X,$ one defines the distance from $x$ to $A$ as
$d(x,A):=\inf\{\left\Vert x-a\right\Vert \mid a\in A\}.$ As usual, we use the
convention $d(x,\emptyset)=\infty.$ The distance function to $A$ is defined as
$d_{A}=d(\cdot,A).$ For a non-empty set $A\subset X,$ we put
$\operatorname*{int}A$ for the topological interior. When we work on a product
space, we consider the sum norm, unless otherwise stated.

Consider now a multifunction $F:X\rightrightarrows Y$ between the Banach
spaces $X$ and $Y$. The domain and the graph of $F$ are denoted respectively
by
\[
\operatorname*{Dom}F:=\{x\in X\mid F(x)\neq\emptyset\}
\]
and%
\[
\operatorname*{Gr}F=\{(x,y)\in X\times Y\mid y\in F(x)\}.
\]
If $A\subset X$ then $F(A):=%
{\displaystyle\bigcup\limits_{x\in A}}
F(x).$ The inverse set-valued map of $F$ is $F^{-1}:Y\rightrightarrows X$
given by $F^{-1}(y)=\{x\in X\mid y\in F(x)\}$. Recall that $F$ is said to be
inner semicontinuous at $(\overline{x},\overline{y})\in\operatorname*{Gr}F$ if
for every open set $D\subset Y$ with $\overline{y}\in D,$ there exists a
neighborhood $U\in\mathcal{V}(\overline{x})$ such that for every $x\in U,$
$F(x)\cap D\neq\emptyset$ (where $\mathcal{V}(\overline{x})$ stands for the
system of the neighborhoods of $\overline{x}$). On the other hand, $F$ is said
to be Lipschitz-like around $(\overline{x},\overline{y})$ with constant $L>0$
if there exist two neighborhoods $U\in\mathcal{V}(\overline{x}),$
$V\in\mathcal{V}(\overline{y})$ such that, for every $x,u\in U,$%
\begin{equation}
F(x)\cap V\subset F(u)+L\left\Vert x-u\right\Vert \mathbb{D}_{Y}.
\label{LLip_like}%
\end{equation}

\bigskip

Let us begin our study with the presentation of the scalar case. Our
motivation is provided by some parametric optimization problems (in our
notations the set of parameters is $P$ and it is initially taken as a
topological space). In the simplest case of a scalar objective defined by a
parametric function $f:X\times P\rightarrow\mathbb{R}$, we look at the problem%
\[
(P_{S}):\text{ }\min f(x,p)\text{ subject to }0\in H(x,p),
\]
where $H:X\times P\rightrightarrows Y$ is a multifunction which defines a
generalized constraints system by the relation $0\in H(x,p).$ One can find in
literature this kind of problems under the generic term of "optimization with
equilibrium constraints". The implicit set-valued map $S:P\rightrightarrows X$
associated to $H$ is%
\[
S(p)=\{x\in X\mid0\in H(x,p)\}.
\]
In fact, if for any fixed $p\in P$ one considers the problem of minimizing the
function $f(\cdot,p)$ with the constraint $0\in H(x,p),$ the set $S(p)$ is the
feasible set of this problem. In general, the treatment of the problem would
involve calculus associated to the set-valued map $S,$ but this becomes a
rather delicate situation since, in general, $S$ and the associated
coderivatives could be hard to compute. One possibility would be as in
\cite[Theorem 5.3]{DurStr2}, and this involves implicit multifunction
theorems, in the line of those presented in Section 3 of this paper.

Let us consider first the case of parametric paradigm, i.e. the case where one
considers the optimality with respect to $x$ and for some fixed values of
parameters. In this respect, let us take $M\subset P$ as a nonempty set. We
say that $\overline{x}\in X$ is a solution for $(P_{S})$ with respect to $M$
if $\overline{x}\in%
{\displaystyle\bigcap\nolimits_{p\in M}}
S(p)$ and for every $p\in M$ there exists $\varepsilon_{p}>0$ s.t. for every
$x\in B(\overline{x},\varepsilon_{p})\cap S(p),$ one has
\begin{equation}
f(\overline{x},p)\leq f(x,p). \label{rminu}%
\end{equation}
This definition actually says that for every fixed $\widetilde{p}\in M,$
$\overline{x}$ is a local solution for the scalar non-parametric problem
\[
\min f(x,\widetilde{p})\text{ subject to }0\in H(x,\widetilde{p}).
\]

In order to illustrate this definition, let us consider $f:\mathbb{R\times
R\rightarrow R},$ $f(x,p)=x^{2}+p^{2}$ and $H:\mathbb{R\times
R\rightrightarrows R},$ $H(x,p)=-x-p+1+\mathbb{R}_{+}.$ Now, $\overline{x}=0$
is a solution for this $(P_{S})$ with respect to $M:=[-1,+\infty).$

Now, we use the Clarke penalization technique for Lipschitz functions (see
\cite[Proposition 2.4.3]{Cla}) in our context. For this, we need to define a
concept of equi-lipschitzianity for the parametric objective function. In the
above notation, for a positive $L,$ one says that the function $f$ is
$L-$Lipschitz at $\overline{x}\in X$ with respect to $M\subset P$ if for every
$p\in M$ there exists a neighborhood $U_{p}\in\mathcal{V}(\overline{x}),$ such
that, for every $x,u\in U_{p},$%
\begin{equation}
\left\vert f(x,p)-f(u,p)\right\vert \leq L\left\Vert x-u\right\Vert .
\label{rlipu}%
\end{equation}

The proof of the next result is given only for completeness.

\begin{thm}
\label{tCl1}Suppose that $f$ is $L$-Lipschitz at $\overline{x}\in X$ with
respect to $M\subset P$ and $\overline{x}$ is solution for $(P_{s})$ with
respect to $M$. Then for every $p\in M,$ there exists a neighborhood $V_{p}$
of $\overline{x}$ s.t. for every $x\in V_{p},$%
\begin{equation}
f(\overline{x},p)\leq f(x,p)+Ld(x,S(p)). \label{rpen}%
\end{equation}

\end{thm}

\noindent\textbf{Proof.} Let $p\in M.$ Let $U_{p}$ be a neighborhood of
$\overline{x}$ s.t. both relations (\ref{rminu}) and (\ref{rlipu}) hold. Now
the arguments follows as in \cite[Proposition 2.4.3]{Cla}. There exists
$\theta>0$ s.t. $B(\overline{x},\theta)\subset U_{p}.$ Consider $V_{p}%
=B(\overline{x},\theta/3)$ and take $x\in V_{p}.$ Clearly, if $x\in V_{p}$ and
$p\in M$ with $x\in S(p)$, then (\ref{rpen}) holds from the definition of the
solution concept. Consider the situation where $x\in V_{p},p\in M,$ but
$(x,p)\notin\operatorname*{Gr}H.$ Then for every $\varepsilon\in(0,\theta/3)$
there is $x_{\varepsilon}^{p}\in S(p)$ s.t.
\begin{align*}
\left\Vert x-x_{\varepsilon}^{p}\right\Vert  &  <d(x,S(p))+\varepsilon\\
&  \leq\left\Vert x-\overline{x}\right\Vert +\varepsilon\\
&  \leq\theta/3+\varepsilon<2\theta/3.
\end{align*}
Hence,
\begin{align*}
\left\Vert x_{\varepsilon}^{p}-\overline{x}\right\Vert  &  \leq\left\Vert
x_{\varepsilon}^{p}-x\right\Vert +\left\Vert x-\overline{x}\right\Vert \\
&  <2\theta/3+\theta/3=\theta.
\end{align*}
Consequently, $x_{\varepsilon}^{p}\in U_{p}\cap S(p),$ so,
\begin{align*}
f(\overline{x},p)  &  \leq f(x_{\varepsilon}^{p},p)\leq f(x,p)+L\left\Vert
x-x_{\varepsilon}^{p}\right\Vert \\
&  \leq f(x,p)+L(d(x,S(p))+\varepsilon)\\
&  =f(x,p)+Ld(x,S(p))+L\varepsilon.
\end{align*}
Letting $\varepsilon\rightarrow0,$ we obtain the conclusion.$\hfill\square$

\medskip

In scalar non-parametric constraint optimization this kind of penalization
with a distance function is very useful when one additionally uses a metric
regularity property of the constraint system. In our generalized setting the
needed regularity is to exist $r>0$ s.t. an inequality of the form%
\begin{equation}
d(x,S(p))\leq rd(0,H(x,p)) \label{metreg}%
\end{equation}
holds for every $p\in M$ and for all $x$ in a neighborhood of $\overline{x}.$
As we shall see in the fourth section, such a relation allows us to work with
the initial set-valued map $H$ instead of the implicit set-valued map $S.$ The
study of this relation is also one the main topics of the present paper and we
develop it in the third section. Later on, we shall be back to the
optimization problems in order to apply the main results.

\bigskip

Moreover, let us consider another possibility to define a concept of solution
for $(P_{S}),$ corresponding to the case where we think $(P_{S})$ as an
optimization problem in both variables, the so-called optimization with
equilibrium constraints. Namely, one says that $(\overline{x},\overline{p})$
is a local solution for $(P_{S})$ if there exists some neighborhoods $U$ and
$W$ of $\overline{x}\ $and $\overline{p},$ respectively, s.t. for all
$(x,p)\in U\times W$ with $0\in H(x,p),$ one has:
\begin{equation}
f(\overline{x},\overline{p})\leq f(x,p). \label{rminv}%
\end{equation}

The transformation of this problem with constraints into an unconstrained
problem can be done by means of the following result.

\begin{thm}
\label{tCl2}Suppose that $f$ is $L-$Lipschitz at $(\overline{x},\overline
{p})\in X\times P$ and $(\overline{x},\overline{p})$ is a local solution for
$(P_{S})$. Then $(\overline{x},\overline{p})$ is a local minimum of the scalar
function%
\[
(x,p)\rightarrow f(x,p)+Ld((p,x),\operatorname*{Gr}S).
\]

\end{thm}

\noindent\textbf{Proof.} Let $U,W$ be the neighborhoods of $\overline{x}$ and
$\overline{p}$ s.t. both relations (\ref{rminv}) and the Lipschitz property
hold. There exists $\theta>0$ s.t. $B(\overline{x},\theta)\times
B(\overline{p},\theta)\subset U\times W.$ Consider $(x,p)\in B(\overline
{x},\theta/6)\times B(\overline{p},\theta/6).$ If $0\in H(x,p)$ (i.e. $x\in
S(p)$) then one obviously has
\[
f(\overline{x},\overline{p})+Ld((\overline{p},\overline{x}),\operatorname*{Gr}%
S)=f(\overline{x},\overline{p})\leq f(x,p)+Ld((p,x),\operatorname*{Gr}%
S)=f(x,p).
\]
Consider the situation where $(x,p)\in B(\overline{x},\theta/6)\times
B(\overline{p},\theta/6),$ but $(x,p,0)\notin\operatorname*{Gr}H.$ Then for
every $\varepsilon\in(0,\theta/3)$ there is $(p_{\varepsilon},x_{\varepsilon
})\in\operatorname*{Gr}S$ s.t.
\begin{align*}
\left\Vert (x,p)-(x_{\varepsilon},p_{\varepsilon})\right\Vert  &
<d((p,x),\operatorname*{Gr}S)+\varepsilon\\
&  \leq\left\Vert (x,p)-(\overline{x},\overline{p})\right\Vert +\varepsilon\\
&  \leq\theta/3+\varepsilon<2\theta/3.
\end{align*}
Therefore,
\begin{align*}
\left\Vert (x_{\varepsilon},p_{\varepsilon})-(\overline{x},\overline
{p})\right\Vert  &  \leq\left\Vert (x_{\varepsilon},p_{\varepsilon
})-(x,p)\right\Vert +\left\Vert (x,p)-(\overline{x},\overline{p})\right\Vert
\\
&  <2\theta/3+\theta/3=\theta.
\end{align*}
Consequently, $(x_{\varepsilon},p_{\varepsilon})\in\left(  U\times W\right)
\cap\operatorname*{Gr}S^{-1},$ so,
\begin{align*}
f(\overline{x},\overline{p})  &  \leq f(x_{\varepsilon},p_{\varepsilon})\leq
f(x,p)+L\left\Vert (x,p)-(x_{\varepsilon},p_{\varepsilon})\right\Vert \\
&  \leq f(x,p)+L(d((p,x),\operatorname*{Gr}S)+\varepsilon)\\
&  =f(x,p)+Ld((p,x),\operatorname*{Gr}S)+L\varepsilon.
\end{align*}
Letting $\varepsilon\rightarrow0,$ we obtain the conclusion.$\hfill\square$

\medskip

In this case, in order to get reasonable optimality conditions, one needs to
avoid the implicit multifunction $S,$ therefore one imposes graphical
regularity: to exist $r>0$ s.t. an inequality of the form%
\begin{equation}
d((p,x),\operatorname*{Gr}S)\leq rd((x,p,0),\operatorname*{Gr}H) \label{grreg}%
\end{equation}
holds for all $(x,p)$ in a neighborhood of $(\overline{x},\overline{p})$.

This relation is studied, together with (\ref{metreg}), in Section 3. Note
that in \cite{LiuYeZhu} the Clarke penalization is used together with some
other inequalities in order to get optimality conditions.

\section{Concepts and tools}

Most of the results of this paper work for several types of generalized
differentiation objects as we shall made precise later. But, for the clarity
of our discussion we mainly use the constructions developed by Mordukhovich
and his collaborators (see \cite{Mor2006}). We briefly remind these concepts
and results. Firstly, recall that $X^{\ast}$ denotes the topological dual of
the Banach space $X,$ while the symbol $w^{\ast}$ is used for the weak-star
topology of the dual system $(X,X^{\ast}).$

\begin{df}
Let $S$ be a non-empty subset of $X$ and let $x\in S,\varepsilon\geq0.$ The
set of $\varepsilon-$normals to $S$ at $x$ is%
\begin{equation}
\widehat{N}_{\varepsilon}(S,x):=\left\{  x^{\ast}\in X^{\ast}\mid
\underset{u\overset{S}{\rightarrow}x}{\lim\sup}\frac{x^{\ast}(u-x)}{\left\Vert
u-x\right\Vert }\leq\varepsilon\right\}  . \label{eps-no}%
\end{equation}

If $\varepsilon=0,$ the elements in the right-hand side of (\ref{eps-no}) are
called Fr\'{e}chet normals and their collection, denoted by $\widehat
{N}(S,x),$ is the Fr\'{e}chet normal cone to $S$ at $x.$

Let $\overline{x}\in S.$ The basic (or limiting, or Mordukhovich) normal cone
to $S$ at $\overline{x}$ is
\[
N(S,\overline{x}):=\{x^{\ast}\in X^{\ast}\mid\exists\varepsilon_{n}%
\downarrow0,x_{n}\overset{S}{\rightarrow}\overline{x},x_{n}^{\ast}%
\overset{w^{\ast}}{\rightarrow}x^{\ast},x_{n}^{\ast}\in\widehat{N}%
_{\varepsilon_{n}}(S,x_{n}),\forall n\in%
\mathbb{N}
\}.
\]

\end{df}

If $X$ is an Asplund space (i.e. a Banach space where every convex continuous
function is generically Fr\'{e}chet differentiable), the formula for the basic
normal cone takes a simpler form, namely:%
\[
N(S,\overline{x})=\{x^{\ast}\in X^{\ast}\mid\exists x_{n}\overset
{S}{\rightarrow}\overline{x},x_{n}^{\ast}\overset{w^{\ast}}{\rightarrow
}x^{\ast},x_{n}^{\ast}\in\widehat{N}(S,x_{n}),\forall n\in%
\mathbb{N}
\}.
\]

Let $f:X\rightarrow\overline{\mathbb{R}}$ be finite at $\overline{x}\in X;$
the Fr\'{e}chet subdifferential of $f$ at $\overline{x}$ is the set
\[
\widehat{\partial}f(\overline{x}):=\{x^{\ast}\in X^{\ast}\mid(x^{\ast}%
,-1)\in\widehat{N}(\operatorname*{epi}f,(\overline{x},f(\overline{x})))\}
\]
and the basic (or limiting, or Mordukhovich) subdifferential of $f$ at
$\overline{x}$ is%

\[
\partial f(\overline{x}):=\{x^{\ast}\in X^{\ast}\mid(x^{\ast},-1)\in
N(\operatorname*{epi}f,(\overline{x},f(\overline{x})))\},
\]
where $\operatorname*{epi}f$ denotes the epigraph of $f.$ On Asplund spaces
one has
\[
\partial f(\overline{x})=\limsup_{x\overset{f}{\rightarrow}\overline{x}%
}\widehat{\partial}f(x),
\]
and, in particular, $\widehat{\partial}f(\overline{x})\subset\partial
f(\overline{x}).$ If $f$ is convex, then both these subdifferential do
coincide with the classical Fenchel subdifferential. If $\delta_{\Omega}$
denotes the indicator function associated with a nonempty set $\Omega\subset
X$ (i.e. $\delta_{\Omega}(x)=0$ if $x\in\Omega,$ $\delta_{\Omega}(x)=\infty$
if $x\notin\Omega$ ), then for any $\overline{x}\in\Omega,$ $\widehat
{\partial}\delta_{\Omega}(\overline{x})=\widehat{N}(\Omega,\overline{x})$ and
$\partial\delta_{\Omega}(\overline{x})=N(\Omega,\overline{x}).$ Let
$\Omega\subset X$ be a nonempty closed set and take $\overline{x}\in\Omega;$
then one has
\begin{equation}
N(\Omega,\overline{x})=\underset{\lambda>0}{%
{\displaystyle\bigcup}
}\lambda\partial d(\cdot,\Omega)(\overline{x}). \label{distgr}%
\end{equation}

The basic subdifferential satisfies a robust sum rule (see \cite[Theorem
3.36]{Mor2006}): if $X$ is Asplund, $f_{1},f_{2},...,f_{n-1}:X\rightarrow
\mathbb{R}$ are Lipschitz around $\overline{x}$ and $f_{n}:X\rightarrow
\overline{\mathbb{R}}$ is lower semicontinuous around this point, then
\begin{equation}
\partial(\sum_{i=1}^{n}f_{i})(\overline{x})\subset\sum_{i=1}^{n}\partial
f_{i}(\overline{x}). \label{sumrule}%
\end{equation}

We shall also need a calculus rule concerning the partial subgradients, that
one can find in \cite[Corollary 3.44]{Mor2006}. We reproduce here this result
in the (less general) form that we actually need in the sequel.

\begin{pr}
\label{partsub}Let $X,Y$ be Asplund spaces and $\varphi:X\times Y\rightarrow
\overline{\mathbb{R}}=\mathbb{R\cup\{+\infty\}}$ be a Lipschitz function
around $(\overline{x},\overline{y})\in X\times Y.$ Then:%
\[
\partial\varphi(\cdot,\overline{y})(\overline{x})\subset\{x^{\ast}\in X^{\ast
}\mid\exists y^{\ast}\in Y^{\ast}\text{ with }(x^{\ast},y^{\ast})\in
\partial\varphi(\overline{x},\overline{y})\}.
\]

\end{pr}

\begin{df}
Let $F:X\rightrightarrows Y$ be a set-valued map and $(\overline{x}%
,\overline{y})\in\operatorname*{Gr}F.$ Then the Fr\'{e}chet coderivative at
$(\overline{x},\overline{y})$ is the set-valued map $\widehat{D}^{\ast
}F(\overline{x},\overline{y}):Y^{\ast}\rightrightarrows X^{\ast}$ given by
\[
\widehat{D}^{\ast}F(\overline{x},\overline{y})(y^{\ast}):=\{x^{\ast}\in
X^{\ast}\mid(x^{\ast},-y^{\ast})\in\widehat{N}(\operatorname{Gr}%
F,(\overline{x},\overline{y}))\}.
\]

Similarly, the normal coderivative of $F$ at $(\overline{x},\overline{y})$ is
the set-valued map $D_{N}^{\ast}F(\overline{x},\overline{y}):Y^{\ast
}\rightrightarrows X^{\ast}$ given by
\[
D_{N}^{\ast}F(\overline{x},\overline{y})(y^{\ast}):=\{x^{\ast}\in X^{\ast}%
\mid(x^{\ast},-y^{\ast})\in N(\operatorname{Gr}F,(\overline{x},\overline
{y}))\}.
\]

\end{df}

Note that, in fact, the concept of normal coderivative, independently of the
normal cone used in its definition, was introduced in \cite{Mor1980}.

Besides (\ref{distgr}) we shall need the following formula obtained by
Thibault in \cite{Thi} (see also \cite{MorNam} for further generalizations and
details): if $X,Y$ are Banach spaces, $F:X\rightrightarrows Y$ has closed
graph and $(\overline{x},\overline{y})\in\operatorname*{Gr}F$ then%
\begin{equation}
N(\operatorname*{Gr}F,(\overline{x},\overline{y}))=\underset{\lambda>0}{%
{\displaystyle\bigcup}
}\lambda\partial\rho_{F}(\overline{x},\overline{y}), \label{distgen}%
\end{equation}
where $\rho_{F}:X\times Y\rightarrow\mathbb{R\cup\{-\infty\}},\rho
_{F}(x,y)=d(y,F(x)).$

\bigskip

In order to conclude this part we remind a subdifferential chain rule (see
\cite[Corollary 3.43]{Mor2006}). Recall (\cite[Definition 3.25]{Mor2006}) that
a function $f:X\rightarrow Y$ is said to be strictly Lipschitz at
$\overline{x}$ if it is locally Lipschitzian around this point and there
exists a neighborhood $V$ of the origin in $X$ s.t. the sequence $(t_{k}%
^{-1}(f(x_{k}+t_{k}v)-f(x_{k})))_{k\in\mathbb{N}}$ contains a norm convergent
subsequence whenever $v\in V,x_{k}\rightarrow\overline{x},$ $t_{k}%
\downarrow0.$ Suppose that $X,Y$ are Asplund spaces. Let $f:X\rightarrow Y$
and $\varphi:Y\rightarrow\mathbb{R}$ s.t. $f$ is strictly Lipschitz at
$\overline{x}\in X$ and $\varphi$ is Lipschitz around $f(\overline{x});$ then
\begin{equation}
\partial(\varphi\circ f)(\overline{x})\subset%
{\displaystyle\bigcup\limits_{y^{\ast}\in\partial\varphi(f(\overline{x}))}}
\partial(y^{\ast}\circ f)(\overline{x}). \label{chain}%
\end{equation}

\section{Regularity of constraint system}

This section is devoted to survey and then to establish several conditions
ensuring the inequalities (\ref{metreg}) and (\ref{grreg}) with a special
emphasis on the case of epigraphical set-valued maps. We look at two different
ways to guarantee the desired relations: firstly, we impose topological
conditions on $H$ and, secondly, we look after coderivative conditions.

Next, we recall a definition.

\begin{df}
Let $L>0,$ $H:X\times P\rightrightarrows Y$ be a multifunction, $((\overline
{x},\overline{p}),\overline{y})\in\operatorname{Gr}H$ and, for every $p\in P,$
denote $H_{p}(\cdot):=H(\cdot,p).$ Then $H$ is said to be open at linear rate
$L>0,$ or $L-$open, with respect to $x$ uniformly in $p$ around $((\overline
{x},\overline{p}),\overline{y})$ if there exist a positive number
$\varepsilon>0$ and some neighborhoods $U\in\mathcal{V}(\overline{x}),$
$V\in\mathcal{V}(\overline{p}),$ $W\in\mathcal{V}(\overline{y})$ such that,
for every $\rho\in(0,\varepsilon),$ every $p\in V$ and every $(x,y)\in
\operatorname*{Gr}H_{p}\cap\lbrack U\times W],$%
\begin{equation}
B(y,\rho L)\subset H_{p}(B(x,\rho)), \label{pLopen}%
\end{equation}

\end{df}

\bigskip

Now, we remind and comment some existing results in this direction. The first
result emphasizes the link between the partial linear openness of the
multifunction $H$ and the inequalities (\ref{metreg}) and (\ref{grreg}). The
first part of this theorem is proved in \cite{DurStr3} in full Banach spaces
framework (but it works also if $P$ is just a topological space), while the
second part comes on the same lines as in \cite[Theorem 5.2]{DurStr2}.

\begin{thm}
\label{impl}Let $X,Y$ be Banach spaces, $P$ be a topological space, $H:X\times
P\rightrightarrows Y$ be a set-valued map which is inner semicontinuous at
$(\overline{x},\overline{p},0)\in\operatorname{Gr}H$. Suppose that $H$ is open
with linear rate $c>0$ with respect to $x$ uniformly in $p$ around
$(\overline{x},\overline{p},0)$. Then there exist $r_{0}>0$ and $U\in
\mathcal{V}(\overline{p})$ such that, for every $(x,p)\in B(\overline{x}%
,r_{0})\times U,$%
\begin{equation}
d(x,S(p))\leq c^{-1}d(0,H(x,p)). \label{xSp}%
\end{equation}

If, moreover, $P$ is a metric space, then there exist $\overline{r}%
,\overline{t}>0$ such that, for every $(x,p)\in B(\overline{x},\overline
{r})\times B(\overline{p},\overline{t}),$%
\begin{equation}
d((p,x),\operatorname*{Gr}S)\leq(1+c^{-1})d((x,p,0),\operatorname*{Gr}H).
\label{gr_reg}%
\end{equation}

\end{thm}

\bigskip

The desired inequalities follow as well from the coderivative conditions, as
illustrated in the result below proved in \cite{DurStr2}.

\begin{thm}
\label{mreg}Let $X,Y$ be Asplund spaces, $P$ be a topological space and
$H:X\times P\rightrightarrows Y$ be a set-valued map such that $0\in
H(\overline{x},\overline{p})$. Suppose that the following assumptions are satisfied:

(i) there exists $U_{1}\in\mathcal{V}(\overline{p})$ such that, for every
$p\in U_{1},$ $\operatorname*{Gr}H_{p}$ is closed;

(ii) $H$ is inner semicontinuous at $(\overline{x},\overline{p},0);$

(iii) there exist $r,s,c>0$ and $U_{2}\in\mathcal{V}(\overline{p})$ such that,
for every $p\in U_{2}$, every $(x,y)\in\operatorname*{Gr}H_{p}\cap\lbrack
B(\overline{x},r)\times B(\overline{y},s)]$ and every $y^{\ast}\in Y^{\ast
},x^{\ast}\in\widehat{D}^{\ast}H_{p}(x,y)(y^{\ast}),$
\[
c\left\Vert y^{\ast}\right\Vert \leq\left\Vert x^{\ast}\right\Vert .
\]

Then the following are true:

(a) For every $a\in(0,c),$ there exist $U\in\mathcal{V}(\overline{p})$ and
$\tau>0$ such that, for every $(x,p)\in B(\overline{x},\tau)\times U,$
\begin{equation}
d(x,S(p))\leq a^{-1}d(0,H(x,p)). \label{mrg}%
\end{equation}

(b) If, moreover, $P$ is a metric space, there exist $\gamma_{0}>0,\tau_{0}>0$
such that, for every $(x,p)\in B(\overline{x},\tau_{0})\times B(\overline
{p},\gamma_{0}),$%
\begin{equation}
d((p,x),\operatorname*{Gr}S)\leq(1+a^{-1})d((x,p,0),\operatorname*{Gr}H).
\label{grg}%
\end{equation}

\end{thm}

We look now to the special case of epigraphical multifunctions. Remind that,
for a multifunction $G:X\rightrightarrows Y$ and a closed convex proper cone
$C\subset Y,$ the epigraphical multifunction associated with $G$ is $\tilde
{G}:X\rightrightarrows Y$ given by $\tilde{G}(x):=G(x)+C$ for every $x\in X.$
We denote the dual cone of $C$ by $C^{\ast}:=\{y^{\ast}\mid y^{\ast}%
(y)\geq0,\forall y\in C\}.$

(... de ce se ia ..)

The next theorem will be used in the sequel and it presents sufficient
conditions for the linear openness of the epigraphical multifunction in terms
of the Fr\'{e}chet coderivative of the initial multifunction. Note that the
first part is \cite[Theorem 3.6]{DurStr1}, while the second part could be
easily obtained by inspecting the proof of the first conclusion.

\begin{thm}
\label{open_epi}Let $X,Y$ be Asplund spaces, $G:X\rightrightarrows Y$ be a
set-valued map, $C$ be a proper closed convex cone in $Y$ and $(\overline
{x},\overline{y})\in\operatorname*{Gr}G.$ Suppose that the following
assumptions are satisfied:

(i) $\operatorname*{Gr}G$ is locally closed at $(\overline{x},\overline{y});$

(ii) there exist $r,c>0$ s.t. for every $(x,y)\in\operatorname*{Gr}%
G\cap\lbrack B(\overline{x},r)\times B(\overline{y},r)]$ and every $y^{\ast
}\in C^{\ast}\cap\mathbb{S}_{Y^{\ast}},z^{\ast}\in2c\mathbb{B}_{Y^{\ast}%
},x^{\ast}\in\hat{D}^{\ast}G(x,y)(y^{\ast}+z^{\ast}),$%
\[
c\left\Vert y^{\ast}+z^{\ast}\right\Vert \leq\left\Vert x^{\ast}\right\Vert .
\]

Then for every $a\in(0,c),$ there exists $\varepsilon>0$ such that, for every
$\rho\in(0,\varepsilon],$%
\[
B(\overline{y},\rho a)\subset G(B(\overline{x},\rho))+C\cap B(0,(a+1)\rho
)\subset\tilde{G}(B(\overline{x},\rho)).
\]

If $\operatorname*{Gr}G$ is closed, the conclusion is more precise in the
following sense: for every $a\in(0,c),$ there exists $\varepsilon:=\min\left(
\tfrac{1}{2}\left(  \tfrac{c}{c+1}-\tfrac{a}{a+1}\right)  ,\tfrac{r}%
{(a+1)}\right)  >0$ such that, for every $\rho\in(0,\varepsilon],$%
\[
B(\bar{y},\rho a)\subset\tilde{G}(B(\bar{x},\rho)).
\]

\end{thm}

\bigskip

Based on this result, let us now consider in the light of our aim, the case
when $H$ is given as an epigraphical multifunction, i.e. $H(x,p)=F(x,p)+C,$
for every $(x,p)\in X\times P,$ where $C$ is, as above, a proper closed convex
cone in $Y.$ In a similar way to the technique of proving Theorem \ref{mreg},
one can obtain the next theorem, which shows the desired inequalities, which
involve now $S$ (the solution map associated to $H$) and $F.$

\begin{thm}
\label{mreg_epi}Let $X,Y$ be Asplund spaces, $P$ be a topological space and
$F:X\times P\rightrightarrows Y$ be a set-valued map such that $0\in
F(\overline{x},\overline{p})$. Suppose that the following assumptions are satisfied:

(i) there exists $U_{1}\in\mathcal{V}(\overline{p})$ such that, for every
$p\in U_{1},$ $\operatorname*{Gr}F_{p}$ is closed;

(ii) $F$ is inner semicontinuous at $(\overline{x},\overline{p},0);$

(iii) there exist $r,c>0$ and $U_{2}\in\mathcal{V}(\overline{p})$ such that,
for every $p\in U_{2}$, every $(x,y)\in\operatorname*{Gr}F_{p}\cap\lbrack
B(\overline{x},r)\times B(0,r)]$ and every $y^{\ast}\in\mathbb{S}_{Y^{\ast}%
}\cap C^{\ast},$ every $z^{\ast}\in2c\mathbb{B}_{Y^{\ast}},$ and every
$x^{\ast}\in\widehat{D}^{\ast}F_{p}(x,y)(y^{\ast}+z^{\ast}),$
\[
c\left\Vert y^{\ast}+z^{\ast}\right\Vert \leq\left\Vert x^{\ast}\right\Vert .
\]

Then for every $a\in(0,c),$ there exist $U\in\mathcal{V}(\overline{p})$ and
$\tau>0$ such that, for every $(x,p)\in B(\overline{x},\tau)\times U,$
\begin{equation}
d(x,S(p))\leq a^{-1}d(0,F(x,p)). \label{mrgepi}%
\end{equation}

If, moreover, $P$ is a metric space, then for every $a\in(0,c),$ there exist
$\tau_{0}>0,\gamma_{0}>0$ such that, for every $(x,p)\in B(\overline{x}%
,\tau_{0})\times B(\overline{p},\gamma_{0}),$%
\begin{equation}
d((p,x),\operatorname*{Gr}S)\leq(1+a^{-1})d((x,p,0),\operatorname*{Gr}F).
\label{grgepi}%
\end{equation}

\end{thm}

\noindent\textbf{Proof. }Concerning the proof of (\ref{mrgepi}), fix arbitrary
$a\in(0,c)$ and $\rho\in\left(  0,\min\left(  \frac{1}{2}\left(  \frac{c}%
{c+1}-\frac{a}{a+1}\right)  ,\frac{r}{2(a+1)}\right)  \right)  .$ Using the
inner semicontinuity of $F$ at $(\overline{x},\overline{p}),$ we can find
$U_{0}\in\mathcal{V}(\overline{p})$ and $\nu>0$ such that for every $(x,p)\in
B(\overline{x},\nu)\times U_{0},$%
\begin{equation}
F(x,p)\cap B(0,a\rho)\not =\emptyset. \label{stea}%
\end{equation}
Denote $U:=U_{0}\cap U_{1}\cap U_{2},$ $\tau:=\min(\nu,\frac{r}{2})$ and take
$(x,p)\in B(\overline{x},\tau)\times U.$

If $0\in F(x,p)\subset H(x,p),$ then (\ref{mrgepi}) trivially holds. Suppose
that $0\not \in F(x,p)$ and then, for every $\varepsilon>0,$ we can find
$y_{\varepsilon}\in F(x,p)$ such that%
\begin{equation}
\left\Vert y_{\varepsilon}\right\Vert <d(0,F(x,p))+\varepsilon. \label{stea2}%
\end{equation}

Because from (\ref{stea}) we have that $d(0,F(x,p))<a\rho,$ we can take
$\varepsilon>0$ sufficiently small such that $d(0,F(x,p))+\varepsilon<a\rho.$
Using (\ref{stea2}), we have that%
\[
0\in B(y_{\varepsilon},d(0,F(x,p))+\varepsilon)\subset B(y_{\varepsilon}%
,a\rho).
\]

Moreover,%
\begin{align*}
B(x,2^{-1}r)  &  \subset B(\overline{x},r),\\
B(y_{\varepsilon},2^{-1}r)  &  \subset B(0,2^{-1}r+a\rho)\subset B(0,r).
\end{align*}

Hence we can apply Theorem \ref{open_epi} for $(x,y_{\varepsilon}%
)\in\operatorname*{Gr}F_{p},r_{0}:=2^{-1}r,$ and $\rho_{0}:=\frac{1}%
{a}(d(0,F(x,p))+\varepsilon)<\rho,$ showing that%
\[
B(y_{\varepsilon},a\rho_{0})\subset F_{p}(B(x,\rho_{0}))+C=H_{p}(B(x,\rho
_{0})).
\]

We can find then $\widetilde{x}\in B(x,\rho_{0})$ such that $y\in
H_{p}(\widetilde{x}),$ or $\widetilde{x}\in S(p).$ Hence%
\[
d(x,S(p))\leq\left\Vert x-\widetilde{x}\right\Vert <\tfrac{1}{a}%
(d(y,F(x,p))+\varepsilon).
\]

Making $\varepsilon\rightarrow0,$ we obtain (\ref{mrgepi}).

Let us now prove (\ref{grgepi}). Take as above $a\in(0,c)$ and $\rho\in\left(
0,\min\left(  \frac{1}{2}\left(  \frac{c}{c+1}-\frac{a}{a+1}\right)  ,\frac
{r}{4(a+1)}\right)  \right)  ,$ use again the inner semicontinuity of $F$ at
$(\overline{x},\overline{p})$ and find the neighborhood $U_{0}$ of
$\overline{p}$ and $\nu>0$ such that for every $(x,p)\in B(\overline{x}%
,\nu)\times U_{0},$ (\ref{stea}) holds. If $P$ is a metric space, we can find
$\gamma>0$ such that $B(\overline{p},\gamma)\subset U_{0}\cap U_{1}\cap
U_{2}.$ Take $\tau_{0}:=\min(\frac{\gamma}{3},\nu,\frac{r}{4}),\gamma
_{0}:=\frac{\gamma}{3}$ and choose $(x,p)\in B(\overline{x},\tau_{0})\times
B(\overline{p},\gamma_{0}).$

We have that%
\begin{align*}
d((x,p,0),\operatorname*{Gr}F)  &  \leq\left\Vert x-\overline{x}\right\Vert
+d(p,\overline{p})<\tfrac{\gamma}{3}+\tfrac{\gamma}{3}=\tfrac{2\gamma}{3},\\
d((x,p,0),\operatorname*{Gr}F)  &  <d(0,F(x,p))<a\rho.
\end{align*}

Without loss of generality suppose that $0\not \in F(x,p),$ hence for every
$\varepsilon>0$ sufficiently small such that $d((x,p,0),\operatorname*{Gr}%
F)+\varepsilon<\min(a\rho,\frac{2\gamma}{3})$ we can find $(x_{\varepsilon
},p_{\varepsilon},y_{\varepsilon})\in\operatorname*{Gr}F$ satisfying%
\begin{align}
\max(\left\Vert y_{\varepsilon}\right\Vert ,d(p_{\varepsilon},p))  &
\leq\left\Vert y_{\varepsilon}\right\Vert +\left\Vert x_{\varepsilon
}-x\right\Vert +d(p_{\varepsilon},p)\nonumber\\
&  <d((x,p,0),\operatorname{Gr}F)+\varepsilon. \label{stea3}%
\end{align}

Hence,
\begin{align}
p_{\varepsilon}  &  \in B(p,\tfrac{2\gamma}{3})\subset B(\overline{p}%
,\gamma),\nonumber\\
0  &  \in B(y_{\varepsilon},d((x,p,0),\operatorname{Gr}F)+\varepsilon)\subset
B(y_{\varepsilon},a\rho) \label{stea4}%
\end{align}

\noindent and
\begin{align*}
B(x_{\varepsilon},4^{-1}r)  &  \subset B(x,4^{-1}r+a\rho)\subset
B(x,2^{-1}r)\subset B(\overline{x},r),\\
B(y_{\varepsilon},2^{-1}r)  &  \subset B(0,2^{-1}r+a\rho)\subset B(0,r).
\end{align*}

Then we can apply Theorem \ref{open_epi} for $(x_{\varepsilon},p_{\varepsilon
},y_{\varepsilon})$ such that $y_{\varepsilon}\in F_{p_{\varepsilon}%
}(x_{\varepsilon}),$ $r^{\prime}:=4^{-1}r$ and $\rho^{\prime}:=\frac{1}%
{a}(d((x,p,0),\operatorname{Gr}F)+\varepsilon)<\rho$ and, using (\ref{stea4}),
we obtain that%
\[
0\in B(y_{\varepsilon},a\rho^{\prime})\subset F_{p_{\varepsilon}%
}(B(x_{\varepsilon},\rho^{\prime}))+C\cap B(0,(a+1)\rho)\subset
H_{p_{\varepsilon}}(B(x_{\varepsilon},\rho^{\prime})).
\]

Then we have that there exists $\widetilde{x}\in B(x_{\varepsilon}%
,\rho^{\prime})$ such that $0\in H_{p_{\varepsilon}}(\widetilde{x}),$ or
$(\widetilde{x},p_{\varepsilon})\in\operatorname{Gr}S.$ Hence, using also
(\ref{stea3}),%
\begin{align*}
d((p,x),\operatorname*{Gr}S)  &  \leq\left\Vert \widetilde{x}-x\right\Vert
+d(p_{\varepsilon},p)\\
&  \leq\left\Vert \widetilde{x}-x_{\varepsilon}\right\Vert +\left\Vert
x_{\varepsilon}-x\right\Vert +d(p_{\varepsilon},p)\\
&  <\tfrac{1}{a}(d((x,p,0),\operatorname{Gr}F)+\varepsilon
)+d((x,p,0),\operatorname{Gr}F)+\varepsilon.
\end{align*}

Making again $\varepsilon\rightarrow0,$ we obtain (\ref{grgepi}).\hfill
$\square$

\section{Applications}

We come back to our motivational facts exposed in the first section. We have
previously seen, on one hand, how regularity of parametric systems could be
used for transformation of a constraint problem into an unconstrained one by
means of penalization and, on the other hand, how the regularity could be
obtained in different ways by different types of conditions. In this section
we are going to apply the penalization in order to get optimality conditions
for several solution concepts in parametric solid vector optimization. Before
starting, we refer to the recent works \cite{Mor2009} and \cite{BaoGupMor}
where similar problems were considered in some greater generality. Our aim
here is to sample how to get necessary optimality conditions under metrical
and graphical regularity and, in order to keep the accent on metrical
conditions, we restrict the attention to the solid case (i.e. the case where
the ordering cone has nonempty topological interior).

Let $g:X\times P\rightarrow Z$ be a parametric vector valued function taking
values into the Banach space $Z$ ordered by a closed convex pointed cone $K$
with nonempty interior (i.e. $\operatorname*{int}K\neq\emptyset$). As usual,
the order $\leq_{K}$ on $Z$ associated to the cone $K$ is given by the
equivalence $x\leq_{K}y$ if and only if $y-x\in K.$ We recall that if
$R\subset Z$ is a nonempty set, then a point $r\in R$ is called weak minimal
point for $R$ with respect to $K$ if%
\[
(R-r)\cap-\operatorname*{int}K=\emptyset.
\]
Keeping the other notations from previous sections, the problem we propose is
as follows%
\[
(P_{V}):\text{ }\min g(x,p)\text{ subject to }0\in H(x,p),
\]
where "$\min$" has a double meaning based on the notion of weak minimal point,
as follows.

Let $M\subset P$ be a nonempty set. We say that $\overline{x}\in X$ is a weak
solution for $(P_{V})$ with respect to $M$ if $\overline{x}\in%
{\displaystyle\bigcap\nolimits_{p\in M}}
S(p)$ and for every $p\in M$ there is an $\varepsilon_{p}>0$ s.t. for every
$x\in B(\overline{x},\varepsilon)\cap S(p),$ one has
\begin{equation}
g(x,p)-g(\overline{x},p)\notin-\operatorname*{int}K. \label{rsoluniv}%
\end{equation}

\bigskip

The other definition of solution runs as follows: we say that $(\overline
{x},\overline{p})$ is a local weak solution for $(P_{V})$ if there exists some
neighborhoods $U$ and $W$ of $\overline{x}\ $and $\overline{p},$ respectively,
s.t. for all $(x,p)\in U\times W$ with $0\in H(x,p),$ one has:
\begin{equation}
g(x,p)-g(\overline{x},\overline{p})\notin-\operatorname*{int}K. \label{rvect}%
\end{equation}

\bigskip

Once again, we say that $g$ is $L$-Lipschitz at $\overline{x}$ with respect to
$M\subset P$ if relation (\ref{rlipu}) holds for $g$ with the norm instead of
modulus in the left-hand side. The main tool we use in conjunction with the
coderivative calculus is contained in the next results which is proved, even
in a more general setting, in \cite{DurTam}. In this result $\partial$ denotes
the Fenchel subdifferential of a convex function and $\operatorname{bd}(K)$
denotes the topological boundary of $K$.

\begin{thm}
\label{scal}Let $K\subset Z$ be a closed convex cone with nonempty interior.
Then for every $e\in\operatorname{int}K$ the functional $s_{e}:Z\rightarrow
\mathbb{R}$ given by
\begin{equation}
s_{e}(z)=\inf\{\lambda\in\mathbb{R}\mid\lambda e\in z+K\} \label{eq. 1}%
\end{equation}
is continuous, sublinear, strictly-$\operatorname{int}K$-monotone and:

(i) $\partial s_{e}(0)=\{v^{\ast}\in K^{\ast}\mid v^{\ast}(e)=1\}$;

(ii) for every $u\in Z$, $\partial s_{e}(u)\neq\emptyset$ and
\begin{equation}
\partial s_{e}(u)=\{v^{\ast}\in K^{\ast}\mid v^{\ast}(e)=1,v^{\ast}%
(u)=s_{e}(u)\}. \label{subscal}%
\end{equation}

Moreover, $s_{e}$ is $d(e,\operatorname{bd}(K))^{-1}$--Lipschitz and for every
$u\in Z$ and $v^{\ast}\in\partial s_{e}(u)$, $\left\Vert e\right\Vert
^{-1}\leq\left\Vert v^{\ast}\right\Vert \leq d(e,\operatorname{bd}(K))^{-1}$.

If $A\subset Z$ is a nonempty set s.t. $A\cap(-\operatorname*{int}%
K)=\emptyset,$ then $s_{e}(a)\geq0$ for every $a\in A$.
\end{thm}

For $e\in\operatorname{int}K$ we shall denote $d(e,\operatorname{bd}(K))^{-1}$
by $L_{e}$ (the Lipschitz constant for $s_{e}$).

We present now the first result of this section.

\begin{thm}
Suppose that $X,Y,Z,P$ are Asplund spaces. Suppose that $g$ is $L$--strictly
Lipschitz at $\overline{x}\in X$ with respect to $M\subset P$ and
$\overline{x}$ is a weak solution for $(P_{V})$ with respect to $M$. Moreover,
suppose that, for every $p\in M,$ (\ref{metreg}) holds, $H_{p}$ has closed
graph and $H$ is Lipschitz-like around $(\overline{x},p,0).$ Then for every
$e\in\operatorname*{int}K$ and for every $p\in M$ there exist $z^{\ast}\in
K^{\ast},z^{\ast}(e)=1$ and $y^{\ast}\in Y^{\ast}$ s.t.%
\[
0\in\partial(z^{\ast}\circ g(\cdot,p))(\overline{x})+D^{\ast}H_{p}%
(\overline{x},0)(y^{\ast}).
\]

\end{thm}

\noindent\textbf{Proof.} Take $e\in\operatorname*{int}K$ and $p\in M.$
Firstly, observe that the scalar application $(x,p)\mapsto s_{e}%
(g(x,p)-g(\overline{x},p))$ is $L\cdot L_{e}$-Lipschitz at $\overline{x}\in X$
with respect to $M\subset P$ because for any $u$ in an appropriate
neighborhood of $\overline{x}$ one has
\begin{align*}
\left\vert s_{e}(g(x,p)-g(\overline{x},p))-s_{e}(g(u,p)-g(\overline
{x},p))\right\vert  &  \leq L_{e}\left\Vert g(x,p)-g(u,p)\right\Vert \\
&  \leq LL_{e}\left\Vert x-u\right\Vert .
\end{align*}
Moreover, since $\overline{x}$ is a weak solution for $(P_{V})$ with respect
to $M$ and taking into account the last conclusion of Theorem \ref{scal} one
deduces that $\overline{x}$ is a solution for
\[
\min s_{e}(g(x,p)-g(\overline{x},p))\text{ subject to }0\in H(x,p),
\]
with respect to $M.$ Note that, for any $p\in M,$ the value of this problem is
$0.$ Whence, following Theorem \ref{tCl1}, there exists a neighborhood $V_{p}$
of $\overline{x}$ s.t. for every $x\in V_{p},$%
\[
0\leq s_{e}(g(x,p)-g(\overline{x},p))+LL_{e}d(x,S(p)).
\]
Now, eventually taking a smaller neighborhood of $\overline{x}$ (denoted
$V_{p}$ as well) and using (\ref{metreg}) one gets that for every $x\in
V_{p},$
\[
0\leq s_{e}(g(x,p)-g(\overline{x},p))+LL_{e}rd(0,H(x,p)).
\]
Taking into account that for $x=\overline{x}$ the right-hand side of the above
problem is $0,$ one deduces that $\overline{x}$ is a local minimum of the
scalar function
\[
x\mapsto s_{e}(g(x,p)-g(\overline{x},p))+LL_{e}rd(0,H(x,p)).
\]
Consequently, by the generalized differentiation calculus rules,
\[
0\in\partial\left(  s_{e}(g(\cdot,p)-g(\overline{x},p))+LL_{e}rd(0,H(\cdot
,p))\right)  (\overline{x}),
\]
and because the first function is Lipschitz while the second one is lower
semicontinuous around the reference point $\overline{x}$ (taking into account
that $H$ is Lipschitz-like around $(\overline{x},p,0)$)$,$ we can employ the
calculus rule (\ref{sumrule}) in order to write%
\[
0\in\partial s_{e}(g(\cdot,p)-g(\overline{x},p))(\overline{x})+LL_{e}r\partial
d(0,H(\cdot,p))(\overline{x}).
\]
Since $s_{e}$ is Lipschitz and $g(\cdot,p)$ is strictly Lipschitz, one can
apply (\ref{chain}) to get%

\[
0\in%
{\displaystyle\bigcup\limits_{z^{\ast}\in\partial s_{e}(0)}}
\partial(z^{\ast}\circ g(\cdot,p))(\overline{x})+LL_{e}r\partial
d(0,H(\cdot,p))(\overline{x}),
\]
i.e. there exists $z^{\ast}\in K^{\ast},z^{\ast}(e)=1$ s.t.%
\[
0\in\partial(z^{\ast}\circ g(\cdot,p))(\overline{x})+LL_{e}r\partial
d(0,H(\cdot,p))(\overline{x}),
\]
Note that, from the Lipschitz-like property of $H_{p},$ the scalar function
$d(\cdot,H(\cdot,p))$ is Lipschitz at $(\overline{x},0)$ (see \cite[Section
2]{Roc1985}) and one uses Proposition \ref{partsub} and relation
(\ref{distgr}) to get%
\[
LL_{e}r\partial d(0,H(\cdot,p))(\overline{x})\subset LL_{e}r\{x^{\ast}\in
X^{\ast}\mid\exists u^{\ast}\in Y^{\ast}\text{ with }(x^{\ast},u^{\ast}%
)\in\partial\rho_{H_{p}}(\overline{x},0)\}.\}
\]
Taking into account $(\ref{distgen})$ one concludes the proof.$\hfill\square$

\bigskip

We pass now to the other kind of vector solution we have introduced before.

\begin{thm}
Suppose that $X,Y,Z,P$ are Asplund spaces. Suppose that $g$ is $L$-strictly
Lipschitz at $(\overline{x},\overline{p})$ which is a local weak solution for
$(P_{V}).$ Moreover, suppose that (\ref{grreg}) holds and $H$ has closed
graph. Then for every $e\in\operatorname*{int}K$ there exist $z^{\ast}\in
K^{\ast},z^{\ast}(e)=1,$ $y^{\ast}\in Y^{\ast}$ s.t.%
\[
(0,0)\in\partial(z^{\ast}\circ g(\cdot,\cdot))(\overline{x},\overline
{p})+D^{\ast}H(\overline{x},\overline{p},0)(y^{\ast}).
\]

\end{thm}

\noindent\textbf{Proof.} The proof runs along some similar lines as above. We
point out the main parts. Take $e\in\operatorname*{int}K.$ The scalar
application $(x,p)\mapsto s_{e}(g(x,p)-g(\overline{x},\overline{p}))$ is
$L\cdot L_{e}$-Lipschitz at $(\overline{x},\overline{p}).$ Moreover, one has
that $(\overline{x},\overline{p})$ is a solution for
\[
\min s_{e}(g(x,p)-g(\overline{x},\overline{p}))\text{ subject to }0\in
H(x,p).
\]
Following Theorem \ref{tCl2}, $(\overline{x},\overline{p})$ is a local
solution for the unconstrained problem%
\[
\min\left[  s_{e}(g(x,p)-g(\overline{x},\overline{p}))+LL_{e}%
d((p,x),\operatorname*{Gr}S)\right]  .
\]
Using (\ref{grreg}) one gets that $(\overline{x},\overline{p})$ is a local
solution for (still unconstrained) problem
\[
\min\left[  s_{e}(g(x,p)-g(\overline{x},\overline{p}))+LL_{e}%
rd((x,p,0),\operatorname*{Gr}H)\right]  .
\]
The subdifferential rules applied for the sum of two Lipschitz functions gives%
\[
(0,0)\in\partial s_{e}(g(\cdot,\cdot)-g(\overline{x},\overline{p}%
))(\overline{x},\overline{p})+LL_{e}r\partial d((\cdot,\cdot
,0),\operatorname*{Gr}H)(\overline{x},\overline{p}).
\]
Since $s_{e}$ is Lipschitz and $g(\cdot,p)$ is strictly Lipschitz one can
apply (\ref{chain}). Taking into account (\ref{distgen}) one gets%

\[
(0,0)\in%
{\displaystyle\bigcup\limits_{z^{\ast}\in\partial s_{e}(0)}}
\partial(z^{\ast}\circ g)(\overline{x},\overline{p})+LL_{e}r\partial
d((\cdot,\cdot,0),\operatorname*{Gr}H)(\overline{x},\overline{p}),
\]
i.e. there exists $z^{\ast}\in K^{\ast},z^{\ast}(e)=1$ s.t.%
\[
(0,0)\in\partial(z^{\ast}\circ g)(\overline{x},\overline{p})+LL_{e}r\partial
d((\cdot,\cdot,0),\operatorname*{Gr}H)(\overline{x},\overline{p})\neq
\emptyset.
\]
Now, taking into account that the distance function is $1$-Lipschitz,
following Proposition \ref{partsub},%
\[
\partial d((\cdot,\cdot,0),\operatorname*{Gr}H)(\overline{x},\overline
{p})\subset\{(x^{\ast},p^{\ast})\mid\exists u^{\ast}\in Y^{\ast}\text{ with
}(x^{\ast},p^{\ast},u^{\ast})\in\partial d_{\operatorname*{Gr}H}(\overline
{x},\overline{p},0)\}.
\]
But, from (\ref{distgr}),
\[
\mathbb{R}_{+}\partial d_{\operatorname*{Gr}H}(\overline{x},\overline
{p},0)\subset N(\operatorname*{Gr}H,(\overline{x},\overline{p},0))
\]
whence, replacing $u^{\ast}$ by $-y^{\ast}$
\[
(0,0)\in\partial(z^{\ast}\circ g)(\overline{x},\overline{p})+D^{\ast
}H(\overline{x},\overline{p},0)(y^{\ast}),
\]
and this concludes the proof.$\hfill\square$

\bigskip

\bigskip

Note that in the case where $g$ is taken as a scalar function (see the
framework and notations of Section 1), then one does not need to scalarize the
function and, consequently, the Lipschitz property is enough.\textbf{ }We
would like also to mention that in this situation one does not need all the
calculus rules we quote for Mordukhovich differentiation and, consequently,
this result could be formulated as well for many other subdifferential and
related coderivatives: see \cite{DurTam} for further details and comments.
However, for completeness, we present such a result.

\begin{thm}
Suppose that $X,Y,Z,P$ are Asplund spaces. Suppose that $f$ is $L$--Lipschitz
at $(\overline{x},\overline{p}),$ which is a local weak solution for
$(P_{S}).$ Moreover, suppose that (\ref{grreg}) holds and that $H$ has closed
graph. Then there exists $y^{\ast}\in Y^{\ast}$ s.t.%
\[
(0,0)\in\partial f(\overline{x},\overline{p})+D^{\ast}H(\overline{x}%
,\overline{p},0)(y^{\ast}).
\]

\end{thm}

\bigskip

We conclude by saying that, in view of the schema displayed in this paper,
every result on regularity of the constraint system would bring necessary
optimality conditions for various types of nonsmooth scalar and vector programs.

\end{document}